\documentclass[final, twoside,11pt]{entics}
\usepackage{enticsmacro}
\usepackage{graphicx}
\usepackage[all]{xy}
\usepackage{amsmath}
\usepackage{tikz}
\newcommand{\dda}{\mathord{\mbox{\makebox[0pt][l]{\raisebox{-.4ex}{$\downarrow$}}$\downarrow$}}}
\newcommand{\dua}{\mathord{\mbox{\makebox[0pt][l]{\raisebox{.4ex}{$\uparrow$}}$\uparrow$}}}
\newcommand{\ua}{\mathord{\uparrow}}
\newcommand{\da}{\mathord{\downarrow}}
\newcommand{\rom}[1]{\rm{\uppercase\expandafter{\romannumeral #1}}}

\newcommand{\I}{\item}
\newcommand{\II}{\begin{enumerate}}
\newcommand{\III}{\end{enumerate}}

\def\conf{ISDT 2022} 	%%Fill in the acronym for your conference (with year)
\volume{2}			%Fill in the ENTICS volume number here
\def\volu{2}			% and here
\def\papno{8}			%Fill in your paper number here
\def\mydoi{{\normalshape  \href{https://doi.org/10.46298/entics.10352}{10.46298/entics.10352}}}
\def\copynames{H. Miao, Q. Li, D. Zhao} %% Fill in the first initial and last name of the authors
\def\runtitle{One-step closure, weak one-step closure and meet continuity}

\def\lastname{Miao, Li and Zhao}
\begin{document}
%%%Note the beginning and end of the frontmatter section that starts here%%%%%
\begin{frontmatter}
  \title{One-step Closure, Weak One-step Closure\\ and Meet Continuity}
  %%%%Now the author(s) names(s)%%%%%
  \author{Hualin Miao\thanksref{a}\thanksref{ALL}\thanksref{myemail}}	%%Note NO SPACE between
  \author{Qingguo Li\thanksref{a}\thanksref{ALL}\thanksref{cocemail}}
   \author{Dongsheng Zhao\thanksref{b}\thanksref{coemail}}		%last name and \thanksref{...}
    %%%Next come the addresses%%%%
   \address[a]{School of Mathematics\\ Hunan University\\				%or between \thanksrefs...
    Changsha, Hunan, 410082, China}  							
  \thanks[ALL]{This work is supported by the National Natural Science Foundation of China (No.12231007) and by Hunan Provincial Innovation Foundation For Postgraduate (CX20200419)}   %%Test of \thanks[ALL} here..
   \thanks[myemail]{Email: \href{mailto:myuserid@mydept.myinst.myedu} {\texttt{\normalshape
        miaohualinmiao@163.com}}}
   \thanks[cocemail]{Corresponding author, Email:  \href{mailto:liqingguoli@aliyun.com} {\texttt{\normalshape
        liqingguoli@aliyun.com}}}
   %%%Note: if both authors share same institution, the only list the address once, after the second
   %%%author.
   %%%There also is a link from the first author to the co-author's address to show how to list
   %%%affiliations to more than one institution, when needed.
  \address[b]{School of Mathematics\\Nanyang Technological University\\
    1 Nanyang Walk, Singapore 637616}
  \thanks[coemail]{Email:  \href{mailto:couserid@codept.coinst.coedu} {\texttt{\normalshape
        dongsheng.zhao@nie.edu.sg}}}
\begin{abstract}
   This paper studies the weak one-step closure and one-step closure properties concerning the structure of Scott closures. We deduce that every quasicontinuous domain has weak one-step closure and show that a quasicontinuous poset need not have weak one-step closure. We also constructed  a non-continuous poset with one-step closure, which gives  a negative answer to an open problem posed by Zou et al.. Finally, we investigate the relationship between weak  one-step closure property and  one-step closure property and prove that a poset has one-step closure if and only if it is meet continuous and has weak one-step closure.
\end{abstract}
\end{frontmatter}
\section{Introduction}\label{intro}
The Scott topology is an intrinsic topology on posets, which is the most important topology in domain theory. Scott proved that a domain endowed with the Scott topology is sober. It is well known that a poset is continuous if and only if its Scott closed set lattice is a completely distributive lattice. In \cite{Zhao}, Zhao introduced  the weak one-step closure property  in order to obtain some  characterizations of $Z$-continuous posets. In \cite{Zou}, Zou et al. proposed the one-step closure property  and proved that every  continuous poset has one-step closure. They asked whether all posets with one-step closure are continuous. Since every continuous poset is quasicontinuous, it is natural to wonder whether every quasicontinuous poset also has one-step closure.

In this paper we shall answer the above problems and investigate  other aspects of weak one-step closure and one-step closure properties. We give the outline of this paper below.

In Section 3, we prove that every quasicontinuous domain has weak one-step closure and show, by a counterexample, that a quasicontinuous poset may not have weak one-step closure. In Section 4, we give a negative answer to the problem posed by Zou et al. in \cite{Zou}. In Section 5, we prove that a poset has one-step closure if and only if it is meet continuous and has weak one-step closure.

Some problems are posed for further investigation.
\section{Preliminaries}

 We now recall some basic notions and results to be used later. We refer the readers to \cite{nht2}, \cite{clad3} for more about these.

 Let $P$ be a poset. For any subset $A$ of $P$, let $\ua A=\{y\in P: x\le y \mbox{ for some } x\in A\}$ and $\da A=\{y\in P: y\le x \mbox{ for some } x\in A\}$. A nonempty  subset $D$ of $P$ is \emph{directed}, denoted by $D\subseteq^{\uparrow} P$,  if  every  finite subset of $D$ has an upper bound in $D$.  The \emph{supremum (infimum)} of a subset $A$ of $P$, if exists, means the least upper (greatest lower) bound of $A$ in $P$ and will be denoted by $\sup A$ ($\inf A$, resp.)  A \emph{semilattice} is a poset in which every nonempty finite subset has an inf; the dual notion is the  \emph{sup semilattice}. A \emph{Scott open} subset of $P$ is an upper set $U$ ($U=\ua U$) of $P$ such that, for every directed subset $D$ of $P$ such that $\sup D$ exists and is in $U$, there is a $d\in D$ such that $d \in U$. The complements of  Scott open sets are called \emph{Scott closed sets}. The collection of all Scott open subsets of $P$ form a topology on $P$, which is called the \emph{Scott topology} of $P$ and denoted by $\sigma(P)$. The collection of all Scott closed subsets of $P$ is denoted by $\Gamma (P)$. The space $(P,\sigma(P))$ is simply written as $\Sigma P$. For any $A\subseteq P$, we write $cl(A)$ as the Scott closure of $A$ (the closure of $A$ with respect to the Scott topology). We denote the set of all finite subsets of a poset $P$ by $Fin(P)$. The Smyth preoder on the set of all subsets of $P$ is given by $G \leq H$ if $\ua H \subseteq \ua G$. We say that $G$ is way below $H$ and write $G \ll H$ if for every directed subset $D \subseteq P$, $\sup D \in \ua H$ implies $D\cap  \ua G \neq \emptyset$. We write $G \ll x$ for $G \ll \{x\}$ and $\dua G=\{x\in L\mid G\ll x\}$. For  $x, y \in P$, $x$ is \emph{way-below} $y$, denoted by $x \ll y$, if for any directed subset $D$ of $P$ for which $\sup D$ exists, $y \leq \sup D$ implies $D \cap \ua x \neq \emptyset$. The poset $P$ is \emph{continuous} if  for all $x\in P$, $\dda x=\{y\in L\mid y\ll x\}$ is directed and $x=\sup \dda x$.

 A poset $P$ is \emph{directed complete} if $\sup D$ exists for all $D\subseteq^{\uparrow} P$. A directed complete poset will be called a dcpo.

A subset $A$ of a topological space is saturated if $A$  is the intersection of all open sets containing $A$.  For a topological space $X$, the set of all compact saturated subsets of $X$ is denoted by $Q(X)$. We write $\mathfrak{K} \subseteq_{flt}Q(X)$ represents that $\mathfrak{K}$ is filtered. We denote the set of all open sets of space $X$ by $\mathcal{O} (X)$. On $Q(X)$, we consider the \emph{upper Vietoris topology} generated by the sets $\Box U =\{K\in Q(X)\mid K\subseteq U\}$, where $U\in \mathcal{O}(X)$.

 \begin{definition}(\cite{Zou})\label{one}
  A poset $P$ is said to have \emph{one-step closure} if
  $cl(A)=A^{'}$
  holds for any $A\subseteq P$, where $A^{'}=\{x\in P\mid \exists D\subseteq ^{\uparrow}\da A, x=\sup D\}$.
\end{definition}
\begin{definition}(\cite{clad3})
 A poset $P$ is meet continuous if for any $x\in P$ and any directed set $D$ of $P$ with $\sup D$ existing,  $x \leq \sup D$ implies  $x\in cl(\da D \cap \da x)$.
\end{definition}

\begin{remark}
For a semilattice $L$, one can prove that it is  \emph{meet continuous} if and only if it satisfies $\inf \{x ,\sup D\} = \sup_{d\in D}\inf\{x,d\}$ for any $x \in L$ and any directed set $D \subseteq L$ with $\sup D$ existing.
\end{remark}

\begin{definition}(\cite{clad3})
A poset $P$ is \emph{quasicontinuous}, if for every $x\in P$,
\II
\I[(1)]  $fin(x)=\{F\mid  F \in Fin(P), F \ll x\}$ is a directed family;
\I[(2)]  $\ua x=\bigcap_{F\in fin(x)}\ua F$ for any $x\in P$.

A quasicontinuous dcpo is called a quasicontinuous domain.

\III
\end{definition}

For any quasicontinuous domain $P$, the family $\{\dua F: F\subseteq P \mbox{ is finite}\}$ is a base of the Scott topology on $P$ (\cite{clad3}).

\begin{definition}(\cite{clad3})
  A space $X$ is \emph{well-filtered} if for each filter basis $\mathcal{C}$ of compact saturated sets of $X$ and each open set $U$ with $\bigcap \mathcal{C} \subseteq U$, there is a $K\in \mathcal{C}$ such that $K\subseteq U$.
\end{definition}

\begin{definition}(\cite{R})
   The set $\mathbb{R}$  of all real numbers equipped with  the topology having $\{[x,y)\mid x<y, x, y\in\mathbb{R}\}$ as a base is called the \emph{Sorgenfrey line}, which is  denoted by $\mathbb{R}_{l}$.
\end{definition}
\section{Weak one-step closure}
By \cite{Zou},  every  continuous poset has one-step closure. However, a quasicontinuous poset may not have one-step closure. In this section, we consider a weaker property, called weaker  one-step closure. We prove that every  quasicontinuous domain has the weak one-step closure, but a quasicontinuous poset need not have this property.

\begin{definition}\label{weak-one}
A poset $P$ is said to have the \emph{weak one-step closure} if for any $A\subseteq P$, it holds that $cl(A)=A^{''}$, where $A^{''}=\{x\in P\mid \exists D\subseteq ^{\uparrow}\da A, ~x\leq\sup D\}$
\end{definition}

\begin{remark}
In \cite{Zhao}, Zhao introduced the Definition \ref{weak-one} for an arbitrary set system, and  called it one-step closure. To be consistent with the paper \cite{Zou}, here we  call this  property  weak one-step closure.
\end{remark}

\begin{theorem}\label{g}
  Every  quasicontinuous dcpo has weak one-step closure.
  \begin{proof} It suffices to show that $cl(A)\subseteq A^{''}$ for any subset $A$ of $L$. To this end, let $x\in cl(A)$, $F\in fin(L)$ with $x\in \dua F$. Then $\dua F$ is Scott open as $L$ is quasicontinuous. Hence $\dua F\cap A\neq \emptyset$, which implies that $F\cap \da A\neq \emptyset$. Thus $(F\cap \da A)_{F\in fin(x)}$ is a filtered family (with respective to the Smyth preorder) of nonempty finite subsets of $L$. By Rudin's Lemma (\cite{clad3}), there exists a directed subset $D$ of $\bigcup_{F\in fin(x)}F\cap \da A$ such that $D\cap (F\cap \da A)\neq \emptyset$ for any $F\in fin(x)$. Also, since $L$ is a quasicontinuous domain, $\{\dua F\mid F\in fin(x)\}$ is a neighborhood basis of $x$. This indicates that $x\in cl(D)=\da \sup D$. Note that $D\subseteq \da A$. We conclude that  $x\in A^{''}$. Hence $cl(A)\subseteq A^{''}$.
  \end{proof}
  \end{theorem}

The following example shows  that the converse conclusion of Theorem \ref{g} is not true.
\begin{example}
  Let $L=(\mathbb{N}\times \mathbb{N})\cup \{\top\}$. Define order $\leq$ on $L$ as follows:
\II
\I[(i)] $(m, n)\leq (s, t)$ if and only if $m=s $ and  $n\leq t$;
\I[(ii)]  $x \le \top$ for all $x \in L$.
\III

It is well known that  $L$ is a dcpo and not quasicontinuous. However, we can easily verify that $L$ has weak one-step closure.

Note that this dcpo $L$ does not have one-step closure.

The dcpo  $L$ is illustrated in Figure 1.

\end{example}

\begin{figure*}[h]
\centering
\begin{tikzpicture}[scale=1.0]

\path (-3,1)   node[left] {(1,1)} coordinate (b);
\fill (b) circle (1pt);
\path (-3,2.2)   node[left] {(1,2)} coordinate (c);
\fill (c) circle (1pt);
\path (-3,3.4)   node[left] {(1,3)} coordinate (d);
\fill (d) circle (1pt);
\path (-2,1)  node[left] {(2,1)} coordinate (e);
\fill (e) circle (1pt);
\path (-2,2.2)  node[left] {(2,2)} coordinate (f);
\fill (f) circle (1pt);
\path (-2,3.4)  node[left] {(2,3)} coordinate (g);
\fill (g) circle (1pt);
\path (-1,5)  node[left] {$\top$} coordinate (l);
\fill (l) circle (1pt);
\path (-1,1)  coordinate (n);
\fill (n) circle (0.1pt);
\path (1,1)  coordinate (m);
\fill (m) circle (0.1pt);
\path (-1,2.2)  coordinate (o);
\fill (o) circle (0.1pt);
\path (1,2.2)  coordinate (p);
\fill (p) circle (0.1pt);
\path (-1,3.4)  coordinate (q);
\fill (q) circle (0.1pt);
\path (1,3.4)  coordinate (r);
\fill (r) circle (0.1pt);

\draw   (b) -- (c)  (c) -- (d)    (e) -- (f)   (f) -- (g) ;
\draw[densely dashed] (d) -- (l) (g) -- (l) (n) -- (m) (o) -- (p) (q) -- (r);
\end{tikzpicture}
\\Fig.1. A non-quasicontinuous domain that has weak one-step closure.
\end{figure*}

The following example shows that a quasicontinuous poset may not have weak one-step closure.
\begin{example}
  Let $L=(\mathbb{N}\times (\mathbb{N}\cup\{\omega\}))\cup \mathbb{N}$. We define an order $\leq$ on $L$ as follows:

For any $x, y\in L$,    $x\leq y$ if and only if  one  of the following holds:
\II
 \I[(i)]  $x=(m,n_{1}), y=(m,n_{2}),n_{1}\leq n_{2}$;
\I[(ii)]  $x=(m,n_{1}), y=(m,\omega)$;
\I[(iii)]  $x, y\in \mathbb{N}$ and $x\leq y$ in $\mathbb{N}$;
\I[(iv)]  $x=(m,n),y\in \mathbb{N},y\geq n,m\geq 2$;
\I[(v)]  $x=(1,n_{1}),y=(m_{2},\omega), m_{2}\geq n_{1}$;
\I[(vi)] $x=(m_{1},n),y=(m_{2},n), m_{1}\leq m_{2},m_{1}\geq 2$.
\III
 $L$ can be illustrated in Figure 2. Then $L$ is a quasicontinuous poset, but $L$ does not have weak one-step closure.

 To see this, first note that $(1,\omega)\in cl(\mathbb{N})=L$ and $(1,\omega)\notin \mathbb{N}^{''}$. Hence, $L$ does not have weak one-step closure. It remains to show that $L$ is quasicontinuous.
 \II
\I[ (i)]  For $(1,\omega)$, we have $\{(1,n)\mid n\in \mathbb{N}\}\subseteq ^{\uparrow}\dda (1,\omega )$ and $(1,\omega)=\sup_{n\in \mathbb{N}}(1, n)$.
 \I[(ii)] For each $(1, n)\in L$. Let $F_{n, m}=\{\{(1,n),(2,m)\}\mid  m\in \mathbb{N}\}$. Then $\{F_{n,m}\mid m\in \mathbb{N}\}\subseteq fin(1,n)$ and is a filtered base  with $\ua (1, n)=\bigcap_{n\in\mathbb{N}}\ua F_{n,m}$.
 \I[(iii)] For each $(m, n)\in L$ with $m\in \mathbb{N}$ and $m \ge 2$, we see easily that $(m, n)\ll (m, n)$. In addition, each $(m, \omega)$ with $m\ge 2$ is the supremum of the directed set $\{(m, n): n\in\mathbb{N}\}$ of compact elements.
 \III
 All these together show that $L$ is quasicontinuous.

\end{example}
\begin{figure*}[h]
\centering
\begin{tikzpicture}[scale=1.0]
\path (1,1)   node[left] {(1,1)} coordinate (a);
\fill (a) circle (1pt);
\path (1,1.7)   node[left] {(1,2)} coordinate (b);
\fill (b) circle (1pt);
\path (1,2.4)   node[left] {(1,3)} coordinate (c);
\fill (c) circle (1pt);
\path (2.5,1.7)   node[left] {(2,1)} coordinate (d);
\fill (d) circle (1pt);
\path (2.5,2.4)  node[left] {(2,2)} coordinate (e);
\fill (e) circle (1pt);
\path (4,2.4)  node[left] {(3,1)} coordinate (g);
\fill (g) circle (1pt);
\path (4,3.1)  node[left] {(3,2)} coordinate (h);
\fill (h) circle (1pt);
\path (5.5,3.1)  node[left] {(4,1)} coordinate (j);
\fill (j) circle (1pt);
\path (5.5,3.8)  node[left] {(4,2)} coordinate (k);
\fill (k) circle (1pt);
\path (1,4.5)  node[left] {(1,$\omega$)} coordinate (m);
\fill (m) circle (1pt);
\path (2.5,5.2)  node[left] {(2,$\omega$)} coordinate (n);
\fill (n) circle (1pt);
\path (4,5.9)  node[left] {(3,$\omega$)} coordinate (o);
\fill (o) circle (1pt);
\path (5.5,6.6)  node[left] {(4,$\omega$)} coordinate (p);
\fill (p) circle (1pt);
\path (8.4,4.5)  node[left] {1} coordinate (q);
\fill (q) circle (1pt);
\path (8.4,5.2)  node[left] {2} coordinate (r);
\fill (r) circle (1pt);

\path (5.7,6.7)   coordinate (x);
\fill (x) circle (0.1pt);
\path (7.9,7.6)   coordinate (y);
\fill (y) circle (0.1pt);
\path (8.4,5.5)   coordinate (z);
\fill (z) circle (0.1pt);
\path (8.4,7.7)   coordinate (f);
\fill (f) circle (0.1pt);
\draw  (a) -- (b)  (b) -- (c) (d) -- (e) (g) -- (h) (j) -- (k)  (q) -- (r) (a) -- (n) (b) -- (o) (c) -- (p) (d) -- (g) (g) -- (j) (e) -- (h) (h) -- (k) (n) -- (o) (o) -- (p) ;
\draw[densely dashed] (c)--(m) (e) -- (n)  (h) -- (o) (k) -- (p) (x) -- (y) (z) -- (f) (k) -- (r) (j) -- (q);
\end{tikzpicture}
\\Fig.2. A quasicontinuous poset does not have weak one-step closure.
\end{figure*}
\section{One-step closure}
In \cite{Zou}, Zou, Li and Ho showed that every  continuous poset has one-step closure. They asked whether $L$ is continuous if  it  has one-step closure. We now give a counterexample for their  problem.  We begin with a lemma which is crucial for further study.
\begin{lemma}\label{d}
If $X$ is a well-filtered space and $Q(X)$ endowed with the upper Vietoris topology is first-countable, then $(Q(X),\supseteq )$ has one-step closure.
\begin{proof}
 Let $\mathcal{A}\subseteq Q(X)$ and $K\in cl(\mathcal{A})$. The fact that $Q(X)$ equipped with the upper Vietoris topology is first-countable implies that there exists a countable neighborhood basis $\mathcal{B}_{K}=\{\Box U_{n}\mid n\in \mathbb{N}\}$ of $K$ and $\Box U_{n+1}\subseteq \Box U_{n}$ for any $n\in \mathbb{N}$.

 Claim 1: $\Box U\subseteq \Box V$ implies  $U\subseteq V$ for any $U, V\in \mathcal{O}(X)$.

 Let $x\in U$. Then $\ua x\in \Box U\subseteq \Box V$. In other words $\ua x \subseteq V$. So $U\subseteq V$ holds.

From Theorem 5.8 in \cite{Xu}, we know that the upper Vietoris topology coincides with the Scott topology on $(Q(X),\supseteq)$. It follows that $\Box U_{n}\cap \mathcal{A}\neq \emptyset$ for any $\Box U_{n}\in \mathcal{B}_{K}$ due to the assumption that $K\in cl(\mathcal{A})$. Choose $K_{n}\in \Box U_{n}\cap \mathcal{A}$ for any $n\in \mathbb{N}$. We define $Q_{n}=K\cup \bigcup_{m\geq n}K_{m}$ for any $n\in \mathbb{N}$.

 Claim 2: $Q_{n}\in Q(X)$ for each  $n\in \mathbb{N}$.

  As a union of saturated sets,  $Q_{n}$ is a saturated. It suffices to verify  that $Q_{n}$ is compact. Let $\{W_{i}: i\in I\}$ be a family of open sets of $X$  such that $Q_{n}\subseteq \bigcup_{i\in I}W_{i}$. Then $K\subseteq \bigcup_{i\in I}W_{i}$. As $K$ is compact,  there exists $F_{1}\in Fin(I)$ such that $K\subseteq \bigcup _{i\in F_{1}}W_{i}$. Then, there exists $\Box U_{n_{0}}\in \mathcal{B}_{K}$ such that $K\in \Box U_{n_{0}}\subseteq \Box \bigcup_{i\in F_{1}}W_{i}$. We  consider  the following two cases:

 Case 1.  $n_{0}\leq n$: For any $m\geq n\geq n_{0}$, then $K_{m}\subseteq U_{m}\subseteq U_{n_{0}}$ by Claim 1. Hence, $Q_{n}\subseteq U_{n_{0}}\subseteq \bigcup_{i\in F_{1}}W_{i}$.

 Case 2. $n_{0}> n$: We can obtain that $\bigcup_{m\geq n_{0}}K_{m}\subseteq U_{n_{0}}\subseteq \bigcup_{i\in F_{1}}W_{i}$ by the similar proof to Case 1. Note that $\bigcup_{i=n}^{n_{0}-1}K_{i}\in Q(X)$ and $\bigcup_{i=n}^{n_{0}-1}K_{i}\subseteq \bigcup_{i\in I}W_{i} $. This means that there exists $F_{2}\in Fin(I)$ such that $\bigcup_{i=n}^{n_{0}-1}K_{i}\subseteq \bigcup_{i\in F_{2}}W_{i}$. Therefore, $Q_{n}\subseteq \bigcup_{i\in F_{1}\cup F_{2}}W_{i}$.

 Claim 3: $K=\sup _{n\in \mathbb{N}}Q_{n}=\bigcap_{n\in \mathbb{N}}Q_{n}$.

 It is easy to see that $K\subseteq \bigcap_{n\in \mathbb{N}}Q_{n}$. For the converse, suppose $x\in \bigcap_{n\in \mathbb{N}}Q_{n}$. We claim that $x\in K$. Assume $x\notin K$. This manifests $\da x\cap K=\emptyset$. In other words, $K\subseteq X\backslash \da x$. It follows that there exists $n\in \mathbb{N}$ such that $K\in \Box U_{n}\subseteq \Box X\backslash \da x$. Through Claim 2, we can conclude that $x\in Q_{n}\subseteq U_{n}\subseteq X\backslash \da x$, which contradicts $x\in \da x$.

  Note that $K_{n}\in \mathcal{A}$ for any $n\in \mathbb{N}$ and $Q_{n}\le K_{n}$ (with respect to the reverse inclusion order). So   $(Q_{n})_{n\in \mathbb{N}}$ is a directed subset of $\da \mathcal{A}$ whose supremum equals $K$. Hence, $Q(X)$ has one-step closure.
\end{proof}

\end{lemma}

The conclusion given in the next theorem answers the question from Zou et al..
\begin{theorem}\label{l}
For  the Sorgenfrey line $\mathbb{R}_{l}$,  $Q(\mathbb{R}_{l})$ has one-step closure and  $Q(\mathbb{R}_{l})$ is not continuous.
\begin{proof}
  By Example 5.18 of \cite{Xu}, we know that the poset $Q(\mathbb{R}_{l})$ is not continuous. The space $\mathbb{R}_{l}$ is Hausdorff, thus well-filtered (every Hausdorff space is sober and every sober space is well-filtered). Hence, by Lemma \ref{d}, $Q(\mathbb{R}_{l})$ has one-step closure.
\end{proof}
\end{theorem}
\section{The relationship between weak one-step closure and one-step closure}
In this section, we investigate the relationship between weak one-step closure and one-step closure.

The following lemma and example justify the term "weak one-step closure".
\begin{lemma}\label{a1}
If a poset $P$ has one-step closure, then it has weak one-step closure.
\begin{proof}
It suffices to prove that $cl(A)=A^{''}$ for any subset $A$ of $P$. From the definition of one-step closure, we have $cl(A)=A^{'}$. One sees obviously that $A^{''}\subseteq cl(A)$. Let $x\in cl(A)$. Then $x\in A^{'}$. It follows that there exists $D\subseteq ^{\uparrow}\da A$ such that $x=\sup D$, i.e., $x\in A^{''}$.
\end{proof}
\end{lemma}

The converse of Lemma \ref{a1} is not true.
\begin{example}
  Let $L=\mathbb{N} \cup \{\omega, a\}
  $, where $\mathbb{N}$ denotes all natural numbers. We define an order $\leq$ on $L$ by $x\leq y$ if and only if:
\II
\I[(i)] $x,y \in \mathbb{N}$ and $x\le y$ holds in $\mathbb{N}$, or
\I[(ii)] $x\in L$ and  $y=\omega$.
\III

Then $L$ can be easily illustrated in Figure 3. and $\mathbb{N}^{'}=\mathbb{N}\cup \{\omega\}\subsetneqq \mathbb{N}^{''}=L$.  Thus $L$  does not have one-step closure.  But it is easy to check that $L$ has weak one-step closure.
\end{example}

\begin{figure*}[h]
\centering
\begin{tikzpicture}[scale=1.0]
\path (0,1)   node[left] {$1$} coordinate (a);
\fill (a) circle (1pt);
\path (0,2)   node[right] {2} coordinate (b);
\fill (b) circle (1pt);
\path (0,3)   node[right] {3} coordinate (c);
\fill (c) circle (1pt);
\path (1,4)   node[right] {a} coordinate (d);
\fill (d) circle (1pt);
\path (0,5)  node[right] {$\omega$} coordinate (e);
\fill (e) circle (1pt);
\draw  (a) -- (b)  (b) -- (c) (d) -- (e);
\draw[densely dashed] (c)--(e);
\end{tikzpicture}
\\Fig. 3. $\mathbb{N}^{''}\neq \mathbb{N}^{'}$.
\end{figure*}

It is then natural to wonder under what conditions, a poset has one-step closure if it has weak one-step closure.
We shall prove that if a poset has weak one-step closure, then it has one-step closure if and only if it is meet continuous.

\begin{lemma}\label{directed}
   Let $L$ be a poset. Then the following statements are equivalent:
\II
\I[(1)] $A^{'}$ is a lower set for any $A\subseteq L$;
\I[(2)]  $D^{'}$ is a lower set for any directed subset $D$ of $L$.
\III

\begin{proof}
     $(1)\Rightarrow (2)$ is trivial.

     $(2)\Rightarrow (1)$ Assume $x\leq y\in A^{'}$. Then there exists $D\subseteq ^{\uparrow } \da A$ such that $y=\sup D$. This means that $x\leq y\in D^{'}$. It follows that $x\in \da D^{'}=D^{'}$. So we have that there exists a directed subset $E$ of $\da D$ such that $x=\sup E$. Note that $E\subseteq \da D\subseteq \da A$. Therefore, $x\in A^{'}$.
   \end{proof}
 \end{lemma}

 If $L$ has one-step closure, then for any subset $A\subseteq L$, $cl(A)=A^{'}$, so it is a lower set.

 In \cite{Zou}, Zou, Li and Ho proved that $L$ is meet continuous if $L$ has one-step closure. We now deduce this result using a weak assumption.

\begin{lemma}\label{D}
 Let $L$ be a poset. If $D^{'}=\da D^{'}$ for any $D\subseteq^{\uparrow} L$, then $L$ is meet-continuous.
 \begin{proof}
  Let $x\in L$, $D\subseteq^{\uparrow} L$ with $\sup D$ existing. If $x\leq \sup D$, then $x\in \da D^{'}=D^{'}$. This means that there exists a directed subset $E$ of $\da D$ such that $x=\sup E$. Note that $E\subseteq \da x\cap \da D$. This implies that $x\in cl(\da x\cap \da D)$. Therefore, $L$ is meet-continuous.
 \end{proof}
\end{lemma}

 \begin{corollary}
 Every poset with one-step closure is meet continuous.
 \end{corollary}

\begin{corollary}\label{inf}
  Let $L$ be a meet continuous semilattice. Then $D^{'}$ is a lower set for any directed subset $D$ of $L$. Moreover, if $L$ has weak one-step closure, then $L$ has one step closure.
\begin{proof}
   From Lemma \ref{lower}, it suffices to prove that $D^{'}$ is a lower set for any directed subset $D$ of $L$. Suppose $x\leq y\in D^{'}$. Then there exists a directed subset $E$ of $\da D$ with $y= \sup E$. The fact that $L$ is a meet continuous semilattice implies that $x=\sup_{e\in E}\inf\{x,e\}$. It is noteworthy that $\{\inf\{x,e\}\mid e\in E\}$ is a directed subset of $\da D$. This means that $x\in D^{'}$. Therefore, $D^{'}$ is a lower set.
\end{proof}
\end{corollary}

 \begin{corollary}\label{sup}
  Let $L$ be a meet continuous sup-semilattice. Then $D^{'}$ is a lower set for any directed subset $D$ of $L$.
  \begin{proof}
  Let $x\leq y \in D^{'}$. This means that there exists a directed subset $E$ of $\da D$ such that $y=\sup E$. Since $L$ is meet continuous,  $x\in cl(\da x\cap \da E)$. It follows that $\da x=cl(\da x\cap \da E)$, and hence $x=\sup (\da x\cap \da E)$.
Let $G=\{\sup F: F\subseteq (\da x\cap \da E) \mbox{ and $F$ is finite}\}$. Then $G$ is a directed subset of  $\da x\cap \da E$ with
$\sup G=\sup (\da x\cap \da E))=x$. Clearly $G\subseteq \da D$, thus $x\in D^{'}$, showing that $D^{'}$ is lower.
  \end{proof}
\end{corollary}

\begin{lemma}\label{lower}
  Let $L$ be a poset with weak one-step closure. If $D^{'}$ is a lower set for any directed subset $D$ of $L$, then $L$ has one-step closure.
  \begin{proof}
  This follows immediately from Definition \ref{weak-one}, Definition \ref{one} and Lemma \ref{directed}.
  \end{proof}
\end{lemma}

\begin{lemma}\label{b1}
  Let $L$ be a meet continuous poset with weak one-step closure. Then $L$ has one-step closure.
\begin{proof}
By Lemma \ref{lower}, it suffices to show that $D^{'}$ is a lower set for any directed $D\subseteq L$.
Suppose $x\leq y\in D^{'}$. Then,  there exists a directed subset $E$ of $\da D$ such that  $y= \sup E$.
Since $L$ is meet continuous, we have  $x\in cl(\da x\cap \da E)$. Since $L$ has weak one-step closure, there is a directed $K\subseteq \da x\cap \da E$ such that $x\le \sup K$. But, trivially $\sup K\le x$, hence $x=\sup K$. In addition, $K\subseteq \da E\subseteq \da D$, so $x\in D^{'}$. Therefore $D^{'}$ is a lower set.
\end{proof}
\end{lemma}

From Lemma \ref{a1} and Lemma \ref{b1} we deduce the following result.
\begin{theorem}
 A  poset has one-step closure if and only if it is meet continuous and  has weak one-step closure.
\end{theorem}

Since every poset having one-step closure is meet continuous, we have the following natural problem.

\begin{problem} Is there a meet continuous poset that does not have one-step closure.
\end{problem}

We have already given in Section 4 an example of a non-continuous poset that has one-step closure. We now confirm that an exact poset with one-step closure is continuous.

\begin{definition}(\cite{Shen}) Let $x, y$ be elements of a poset $P$.
We say that $x$ is \emph{weakly way-below} $y$, denoted by $x \ll_{w} y$, if for any directed subset $D$ of $P$ for which $\sup D$ exists, $y =\sup D$ implies $D \cap \ua x \neq \emptyset$. A poset $P$ is called \emph{exact} if for any $x\in P$, $\dda_{w}x=\{y\in P\mid y\ll_{w} x\}$ is directed and $\sup \dda_{w}x=x$.
\end{definition}
\begin{theorem}
  Let $L$ be a poset. Then the following statements are equivalent:

  $(1)$ $L$ is continuous;

  $(2)$ $L$ has one-step closure and is exact;

  $(3)$ $A^{'}$ is a lower set for any $A\subseteq L$ and $L$ is an exact poset.

  $(4)$ $D^{'}$ is a lower set for any directed subset $D$ of $L$ and $L$ is an exact poset.

  \begin{proof}
    $(1)\Rightarrow(2) \Rightarrow (3) \Rightarrow (4)$ are all obvious.

    $(4)\Rightarrow (1)$ From the definition of exact posets and continuous posets, it suffices to prove that $x\ll_{w}y$ implies that $x\ll y$ for any $x,y\in L$. Let $D$ be a directed subset of $L$ with $y\leq \sup D$, it follows that $\sup D\in D^{'}$. We have that $y\in D^{'}$ since $D^{'}$ is a lower set. This means that there exists $E\subseteq ^{\uparrow}\da D$ such that $y=\sup E$. The assumption that $x\ll_{w}y$ reveals that $\ua x\cap E\neq \emptyset$. Hence, $\ua x\cap D\neq \emptyset$.
  \end{proof}
\end{theorem}

The following is a problem concerning the connections among the concepts of meet continuity, exactness and continuity of posets.

\begin{problem}
  Let $L$ be a meet continuous and exact poset. Must  $L$ be continuous?
\end{problem}
Although we cannot solve the above problems, we have the following corollary by Corollary \ref{inf} and Corollary \ref{sup}.

\begin{corollary}
  Let $L$ be a meet continuous semilattice or sup-semilattice. Then $L$ is continuous iff $L$ is an exact poset.
\end{corollary}

\begin{proposition}\label{k}
  Let $L, M$ be two posets. If $L$ is a Scott retract of $M$, which has one-step closure, then $L$ has one-step closure.
  \begin{proof}
    Since $L$ is a Scott retract of $M$, we have that there exists two Scott continuous maps $s: L\rightarrow M$ and $r: M\rightarrow L$ such that $id_{L}=r\circ s$. Now let $x\in L$, $A\subseteq L$ with $x\in cl(A)$. It follows that $s(x)\in cl (s(A))$ by the Scott continuity of $s$. We know that there exists $D\subseteq ^{\uparrow}\da s(A)$ such that $s(x)=\sup D$ since $M$ has one-step closure. This implies that $x=r\circ s(x)=r(\sup D)=\sup r(D)$ from the Scott continuity of $r$. Note that $r(D)\subseteq \da r(s(A))=\da A$. This means that $L$ has one-step closure.
  \end{proof}
\end{proposition}
\begin{lemma}
 Let $L,M$ be two posets. If $L$ is a Scott retract of $M$, which has weak one-step closure, then $L$ has weak one-step closure.
 \begin{proof}
   The proof is similar to Proposition \ref{k}.
 \end{proof}
\end{lemma}

\bibliographystyle{./entics}

\end{document}